\documentclass[10pt,a4paper]{article}
\usepackage{amssymb}
\usepackage{amsmath}
\usepackage{latexsym}
\usepackage{graphicx}
\usepackage{amsfonts}
\usepackage{epstopdf}
\usepackage{color}
\usepackage{caption}
\usepackage{subcaption}
\usepackage{helvet}
\newtheorem{lemma}{Lemma}
\newtheorem{theorem}{Theorem}

\catcode`\@=11
\catcode`\@=12

\def\E{\mathbb{E\,}}

\newcommand{\Var}{\text{\textnormal{Var}}}
\newcommand{\Cov}{\text{\textnormal{Cov}}}

\makeatletter
\renewcommand\section{\@startsection{section}{1}{\z@}%
                                  {-3.5ex \@plus -1ex \@minus -.2ex}%
                                  {2.3ex \@plus.2ex}%
                                  {\normalfont\large\bfseries}}
\makeatother

\newcommand\dist{\mathrel{\overset{\makebox[0pt]{\mbox{\normalfont\tiny\sffamily d}}}{\approx}}}

\newcommand{\lo}{\text{\textnormal{o}}}
\newcommand{\I}{\text{\textnormal{I}}}
\newcommand{\halmos}{\vspace{3mm} \hfill \mbox{$\Box$}}
\newcommand{\di}{\text{\textnormal{d}}}
\newcommand{\idef}{\stackrel{\mathrm{def}}{=}}
\renewcommand{\tilde}{\widetilde}
\renewcommand{\hat}{\widehat}
\renewcommand{\bar}{\overline}

\begin{document}

\title{A Correction Term for the Covariance of\\
Renewal-Reward Processes with Multivariate Rewards}
\author{
Brendan Patch\thanks{Corresponding author. Email: b.patch@uq.edu.au},\thanks{The School of Mathematics and Physics, The University of Queensland, Qld 4169, Australia.}
\quad Yoni Nazarathy,\footnotemark[2]
\quad and\quad Thomas Taimre.\footnotemark[2]
}


\maketitle

\thispagestyle{empty}

\vspace{-7mm}
\begin{abstract}
 We consider a renewal-reward process with multivariate rewards.
 Such a process is constructed from an i.i.d.\ sequence of time periods, to each of which there is associated a multivariate reward vector.
 The rewards in each time period may depend on each other and on the period length, but not on the other time periods.
 Rewards are accumulated to form a vector valued process that exhibits jumps in all coordinates simultaneously, only at renewal epochs.

We derive an asymptotically exact expression for the covariance function (over time) of the rewards, which is used to refine a central limit theorem for the vector of rewards.
 As illustrated by a numerical example, this refinement can yield improved accuracy, especially for moderate time-horizons.

\end{abstract}

\section{Introduction}

Probabilistic modeling and analysis has a long tradition in dealing with the behaviour of regenerative processes. Such processes restart probabilistically at renewal instances, forming a sequence of
independent and identically distributed (i.i.d.) sub-processes. They are used in stochastic simulation, reliability analysis, actuarial studies, queueing theory, and other aspects of applied probability and statistics. An  illustrative example is that of natural disasters. When a natural disaster occurs there are several simultaneous costs (e.g.\ personal property loss and infrastructure damage), which may be distributed across different locations. These costs (rewards) are typically dependent and may also depend on the time elapsed since the previous disaster.  If we assume that the system resets after such an event then the situation is well described by a renewal-reward process with multivariate rewards --- which we will refer to simply as a {\em multivariate renewal-reward process}.

The multivariate renewal-reward process is constructed on a probability space supporting $\{{\mathbf Z}_n\}_{n=0}^\infty$, a sequence of $(L+1)$-dimensional independent random vectors with possibly dependent coordinates. The first coordinate of ${\mathbf Z}_n$, denoted $T_n$, signifies the time between events, which we call {\em renewals}, and is assumed non-negative. The remaining $L$ coordinates, denoted $X_{1,n},\ldots,X_{L,n}$ are the {\em rewards} and are not sign restricted. Assume that $\{{\mathbf Z}_n\}_{n=1}^\infty$ are i.i.d.\ and, as is standard in renewal theory (see e.g.\ \cite{bookAsmussen2003} or \cite{bookGut2009}), ${\mathbf Z}_0$ may follow a different distribution. We refer to the case of $T_0 \equiv 0$ and all $X_{i,0} \equiv 0$ as \emph{ordinary}; otherwise the process is \emph{delayed}. To avoid trivialities assume that $T_n$ and all $X_{i,n}$ are almost surely not zero for $n\ge 1$.

Let $S_n \idef \sum_{i=0}^n T_i$, so that $\{S_n\}_{n=0}^\infty$ are the renewal times. Taking $N(t) \idef \min \{ n \,: \, S_n > t\}$ the multivariate renewal-reward process is $\{{\mathbf R}(t)\,:\,t\ge0\}$, or simply ${\mathbf R}(\cdot)$, where
\begin{align}
\label{eqn:RRmulti}
{\mathbf R}(t) \idef \left[
\sum_{n=0}^{N(t)-1} X_{1,n}\,
,
\ldots
,
\,
\sum_{n=0}^{N(t)-1} X_{L,n}
\right]\,.
\end{align}
We treat the summations in ${\mathbf R}(t)$ as empty for $t<T_0$, since $N(t)=0$ there.
For $L>2$ the $i$-th coordinate of ${\mathbf R}(\cdot)$ is denoted $R_i(\cdot)$. For $L=2$ the coordinates are $R_x(\cdot)$ and $R_y(\cdot)$, and the $n$-th reward vector $[X_{1,n},\,X_{2,n}]$ is written simply as $[X_n,\,Y_n]$, where we represent all vectors as rows. Note that in the ordinary case, $N(0)=1$ and in the delayed case $N(0)=0$.

We focus on the case of moderate or large $t$ and aim to approximate the distribution of ${\mathbf R}(t)$. In the illustrative example of natural disasters, this is the multivariate distribution describing the different types of losses accumulated during the first $[0,t]$ time units.
In some very special cases the distribution of ${\mathbf R}(t)$ admits an explicit form.  If, for example, the coordinates of ${\mathbf Z}_n$ are mutually independent and $T_n$ is exponentially distributed, then $N(\cdot)$ is a Poisson process and ${\mathbf R}(\cdot)$ is a vector of independent compound Poisson processes. In general, however, the distribution of ${\mathbf R}(t)$ is not easily obtainable, in which case asymptotic approximations become particularly appealing.

Under regularity conditions (described in the next section) it is well known that ${\mathbf R}(t)$ obeys a normal central limit theorem (CLT) as $t \to \infty$, where the mean and covariance terms appearing in the CLT are determined by moments of ${\mathbf Z}_1$. Brown and Solomon further established in \cite{brown1975second} that for a renewal-reward process satisfying suitable regularity conditions with univariate rewards ($L=1$),
\[
\E R_x(t) =  a_x\,t + b_x + \lo(1)\,\quad\text{and}\quad
\Var\big(R_x(t) \big) =  c_x\,t + d_x + \lo(1)\,,
\]
where $\lo(1)$ is a function that vanishes as $t \to \infty$. Here, the constants $a_x$ and $c_x$ are determined by moments (including cross moments) of  $(T_1,X_1)$ and the constants $b_x$ and $d_x$ are determined by moments of $(T_0,X_0)$ and $(T_1,X_1)$. Expressions for $a_x$ and $c_x$, as well as a version of $b_x$ with rewards independent of renewals, were found by Smith in \cite{smith1955}. Subsequently, in \cite{brown1975second}, Brown and Solomon extended to find $b_x$ and $d_x$ for the general univariate renewal-reward process.

The main contribution of the current paper is to generalize the result of \cite{brown1975second} to multivariate rewards. We prove that, under regularity conditions,
\begin{align*}
\Cov \big(R_x(t),\, R_y(t) \big) = c_{x,y}\,t + d_{x,y} + \lo(1)\,.
\end{align*}
As before, $c_{x,y}$ depends on the moments of ${\mathbf Z}_1$ and $d_{x,y}$ depends on the moments of both ${\mathbf Z}_1$ and ${\mathbf Z}_0$. Expressions for $c_{x,y}$ appeared in \cite{smith1955}, although without an explicit proof for this form of the covariance curve.  Our expression for $d_{x,y}$ is new and generalizes $d_x$ of \cite{brown1975second}.

The multivariate CLT for ${\mathbf R}(t)$ first appeared in \cite{smith1955}. The CLT uses a covariance matrix with elements $c_{x,y}$ (or $c_x$ on the diagonal).  Our refined asymptotics suggest an improved approximation to ${\mathbf R}(t)$ based on this CLT, our new $d_{x,y}$ term, and the previously known $d_{x}$ term.  We illustrate the usefulness of this improved approximation in an example. A further (minor contribution) of the current paper is in casting Smith's CLT in a modern form.  A related presentation is in Section~7.4 of \cite{bookWhitt2001}, where functional CLTs are given. The case handled there is one dimensional and assumes rewards are independent of renewals.

The remainder of the paper is structured as follows. In Section~\ref{sec:main} we present our main result on correction terms to the covariance curve of the multivariate renewal-reward process as well as the CLT and the improved approximation to ${\mathbf R}(t)$.  Section~\ref{sec:numeric} demonstrates the usefulness of our correction terms through a numerical illustration. Proofs are in Section~\ref{sec:proofs}.  We conclude in Section~\ref{sec:conclusion}.


\section{Main Results}
\label{sec:main}

Our results are stated in terms of moments (and cross moments) of ${\mathbf Z}_0$ and ${\mathbf Z_1}$. It is useful to denote some of the moments of ${\mathbf Z}_1$ as follows: $\mu_i \idef \E T_1^i$, $\lambda_i \idef \E X_{1}^i$, $\alpha_i \idef \E Y_1^i$, $m_{i,j} \idef \E T_1^i\,X_1^j$, $n_{i,j} \idef \E T_1^i\,Y_1^j$, and $p_{i,j,k} \idef \E T_1^i\,X_{1}^j\,Y_1^k$.
Denote the distribution function of $T_0$ by $F_0(\cdot)$ and that of $T_1$ (and subsequent inter-event times) by $F(\cdot)$. We call $F(\cdot)$ \emph{non-lattice} if the corresponding probability measure $\di F(\cdot)$ is not concentrated on a set of the form $\{\delta, 2\delta, \dots\}$. A distribution function is said to have the stronger property of being \emph{spread out} if $F^{(n)}(\cdot)$, the $n$-th convolution of $F(\cdot)$, has a component that is absolutely continuous (e.g.\ \cite[Sec 7.1]{bookAsmussen2003}).

The growth rate $a_x \idef \lambda_1\, /\, \mu_1$ is well known.
In \cite{brown1975second}, it was further established:
\begin{theorem}[Restatement of \cite{brown1975second}, Lemma 1]
\label{thm1}
For $F(\cdot)$ non-lattice and $\mu_2$, $\lambda_1$, and $m_{1,1}$ finite,
\begin{equation}
\label{eq:ERx}
\E R_x(t) = a_x\,t + b_x + \lo(1)\,,
\end{equation}
where $a_x \idef \lambda_1\,/\,\mu_1$ and $b_x \idef \mu_1^{-1}\,\mu_2\,a_x \,/2-\mu_1^{-1}\,m_{1,1} + \E X_0-a_x\,\E T_0$.
\end{theorem}
In order to state our main result, consider an ordinary renewal-reward process where the rewards are distributed as the product $X_1 Y_1$. For ordinary ${\mathbf R}(\cdot)$ denote the $i$-th reward coordinate by $\mathring{R}_i(\cdot)$. In particular,  for $L=2$, we write $\mathring{R}_x(\cdot)$, $\mathring{R}_y(\cdot)$, and $\mathring R_{xy}(\cdot)$ for the two reward coordinates and the associated product reward coordinate. Applying Theorem~\ref{thm1} above, we have that $\E  \mathring R_{xy}(t) \idef \E \sum_{n=1}^{N(t)-1} X_n Y_n$ can be represented as
\begin{equation}
\label{eq:bxyring}
\E \mathring R_{xy}(t) = a_{xy}\,t + \mathring{b}_{xy} + \lo(1)\,,
\quad
\mbox{with}
\quad
a_{xy}= \mu_1^{-1}\,p_{0,1,1}
\quad
\text{and}
\quad
\mathring{b}_{xy} =\mu_1^{-1}a_{xy}\,\mu_2\,/2-\mu_1^{-1}\,p_{1,1,1}\,.
\end{equation}
Our main theorem is a generalization of the key results in \cite{brown1975second}. It utilizes the expressions for $a_{xy},\, \mathring{b}_{xy}$ as well as $a_x$, $a_y = \alpha_1\,/\,\mu_1,$ and the corresponding correction terms in the ordinary case,
\[
\mathring{b}_x = \mu_1^{-1}\,a_x\,\mu_2\,/2-\mu_1^{-1}\,m_{1,1}
\qquad
\mbox{and}
\qquad
\mathring{b}_y = \mu_1^{-1}\,a_y\,\mu_2\,/2-\mu_1^{-1}\,n_{1,1}.
\]
\begin{theorem}
\label{thm2}
For $F(\cdot)$ spread out and $\mu_3$, $\lambda_2$, $\alpha_2$, $m_{1,2}$, $n_{1,2}$, $p_{1,1,1}$, $\E T_0^2$, and $\E X_0\,Y_0$ finite,
\[
\Cov\big(R_{x}(t),\,R_{y}(t)\big) = c_{x,y}\,t + d_{x,y} + \lo(1)\,,
\]
where
\[
c_{x,y} \idef  \mu_1^{-1} \Cov\big(X_1-a_x\,T_1,\, Y_1 - a_y\,T_1 \big) = a_{xy}+a_x\,\mathring{b}_y+a_y\,\mathring b_x\,.
\]
Further,
\begin{equation}
\label{eq:DD}
d_{x,y} \idef \mathring{d}_{x,y}-c_{x,y}\,\E T_0+a_x\,a_y\,\Var(T_0)+\Cov(X_0,\,Y_0)-a_x\,\Cov(T_0,\,Y_0)-a_y\,\Cov(T_0,\,X_0)\,,
\end{equation}
with\vspace{-3mm}
\begin{align*}
\mathring{d}_{x,y} &\idef  \mathring{b}_x \, \mathring{b}_y + \mathring{b}_{xy}+2 \, a_x \,\ell_y+2 \, a_y \, \ell_x\,,
\end{align*}
\vspace{-5mm}
where,
\begin{align}
\label{eq:3543}
\ell_x &\idef  \mu_1^{-3}\,\lambda_1\,\mu_2^2\,/4-\mu_1^{-2}\,\lambda_1\,\mu_3\,/6+\mu_1^{-1}\,m_{2,1}\,/2-\mu_1^{-2}\,\mu_2\,m_{1,1}\,/2\,, \\
\nonumber
\ell_y  &\idef  \mu_1^{-3}\,\alpha_1\,\mu_2^2\,/4-\mu_1^{-2}\,\alpha_1\,\mu_3\,/6+\mu_1^{-1}\,n_{2,1}\,/2-\mu_1^{-2}\,\mu_2\,n_{1,1}\,/2\,.
\end{align}
\end{theorem}
Note that: (i) as shown in Lemma~\ref{lem:BS3} below, the quantity $\ell_x$ (as well its $y$-counterpart) is in fact the integrated $\lo(1)$ term of \eqref{eq:ERx};
(ii) for $y=x$ Theorem~\ref{thm2} reduces to results of \cite{brown1975second} with $\mathring b_{xy}= \mathring b_{xx} = \mu_1^{-2}\,\mu_2\,\lambda_2\,/2-\mu_1^{-1}\,m_{1,2}$; and (iii) for ordinary ${\mathbf R}(\cdot)$ the terms involving ${\mathbf Z}_0$ vanish, implying $d_{x,y} = \mathring{d}_{x,y}$.

For $L$-dimensional ${\mathbf R}(\cdot)$, we define the matrices and vectors:
\[
{\mathbf a} = [\,a_i\,]_{i=1}^L\,,
\quad
{\mathbf b}=[\,b_i\,]_{i=1}^L\,,
\quad
C = \mu_1^{-1} \Cov \big(
\left[
\,\gamma_i\,
\right]_{i=1}^L
\big)\,,
\quad
\text{and}
\quad
D= [\,d_{i,j}\,]_{i,j =1}^n\,.
\]
Here the elements $a_i$, $b_i$, and $d_{i,j}$ are as defined in Theorems~\ref{thm1}~and~\ref{thm2} above, where $x$ and/or $y$ are replaced by some pair $i,j \in \{1,\ldots,L\}$, and $\gamma_i = X_{1,i} - a_{i}\,T_{1}$ for $i=1,\ldots,L$. The vector ${\mathbf a}$ and the covariance matrix $C$ play a role in the CLT which we state now. The vector ${\mathbf b}$ and our (new contribution) matrix $D$ are the {\em correction terms}.
These appear in the refinement that follows.
\begin{theorem}[Originally in \cite{smith1955}]\label{thm:Smith}
If $F(\cdot)$ is spread out, $\E X_{i,1}\,T_1 < \infty$, and $\E X_{i,1}\,X_{j,1} < \infty$ then the sequence (in $t$) of random vectors,
\[
\left[\frac{R_{1}(t) - a_1\,t}{\sqrt{t}}, \ldots, \, \frac{R_{L}(t) - a_L\,t}{\sqrt{t}}\right],
\quad
t > 0\,,
\]
converges in distribution, as $t\to\infty$, to a zero mean normal random vector with covariance matrix $C$, denoted here by ${\sf {\mathbf N}}\left({\mathbf 0},C\right)$.
\end{theorem}

Motivated by Theorems~\ref{thm1}--\ref{thm:Smith}, we suggest the following \emph{refined normal approximation} to the distribution of the multivariate renewal-reward process at time $t$:
\begin{equation}
\label{eq:approx}
{\mathbf R}(t)  ~ \dist ~  {\sf {\mathbf N}}
\Big(
{\mathbf a}\,t + {\mathbf b}
,~
C \,t + D
\Big)\,.
\end{equation}
Note that the matrix $D$ may not be positive definite (PD) --- in other words not a covariance matrix --- whereas $C$ always is.
When $D$ is not PD, it is easy to see that $C\,t + D$ is PD
for all $t$ greater than some $t_0>0$, and is not PD for all $t \le t_0$.
Consequently, we only suggest \eqref{eq:approx} when $C \,t + D$ is PD (i.e.\ $t>t_0$).

\vspace{-5mm}
\section{Numerical Illustration}
\label{sec:numeric}

To illustrate the applicability of our refined normal approximation \eqref{eq:approx} assume that we wish to evaluate
\[
m(t) \idef \E \min\big\{R_x(t),\, R_y(t)\big\}\,.
\]
A simple expression for this expected minimum is generally not available, but by approximating the distribution of $\big[R_x(t),\, R_y(t)\big]$ as normal
using Theorem~\ref{thm:Smith} we obtain a very good approximation for $m(t)$, which is generally improved using our refinement \eqref{eq:approx}.
The expected minimum of a bivariate normal random vector $[W,\,V]$ is \begin{align}
\label{eq:Hunter}
\Phi\left(\frac{\E(V-W)}{\Var\big(W-V\big)}\right)\E W+\Phi\left(\frac{\E(W-V)}{\Var\big(W-V\big)}\right)\E V-\phi\left(\frac{\E(W-V)}{\Var\big(W-V\big)}\right)\Var\big(W-V\big)\,,
\end{align}
where $\Phi$ and $\phi$ are respectively the cdf and pdf of the standard normal distribution (see e.g.\ \cite{Hunter1974}). Thus, using the mean and variance/covariance expansions from Theorems~\ref{thm1}~and~\ref{thm2} we can (for fixed $t$) combine \eqref{eq:Hunter} with \eqref{eq:approx} to get an explicit approximation of $m(t)$, denoted $\tilde{m}(t)$.
For moderate to large $t$, we expect using $C\,t+D$ as the covariance will yield a better approximation of $m(t)$ than only using $C\,t$.
That is, we expect that the matrix $D$, with our newly found covariance refinement term $d_{x,y}$ on the off-diagonals, will improve the approximation.

As a specific numerical example consider
\[
{\bf Z}_n = \big[T_n,\,X_n,\,Y_n\big] = \big[U_{1,n} + U_{4,n},~ U_{2,n} + U_{4,n},~ U_{3,n}+U_{4,n}\big]\,,~n=0,1,2,\ldots,
\]
where $U_{i,n}$ are all independent exponential random variables with unit mean for $i=1,2,4$ and mean $1/2$ for $i=3$. Now, using the expressions of Theorems~\ref{thm1}~and~\ref{thm2} we obtain
\[
{\mathbf a} =  \left[ \begin{array}{cc} 1, & -1\end{array}\right],
\quad
{\mathbf b} =  \left[ \begin{array}{cc} -1,& -8/7\end{array}\right],
\quad
C = \left[ \begin{array}{cc} 1, & 3/8\\ 3/8,&7/16\end{array}\right],
\quad
\text{and}
\quad
D = \left[ \begin{array}{cc} 1/2, & 1/2\\ 1/2,&13/64\end{array}\right]\,.
\]
In this case the refinement to the covariance curve is only applicable for $t > t_0 = (\sqrt{731}-3)/38 \approx 0.63$, since $Ct+D$ is not a PD matrix when $t\leq t_0$.
We estimated the true $m(t)$ by extensive simulation, taking the mean of the minimum over $10^7$ sample paths of ${\mathbf R}(\cdot)$ as the estimate $\hat m(t)$.

\begin{figure}[h!]
	\centering
	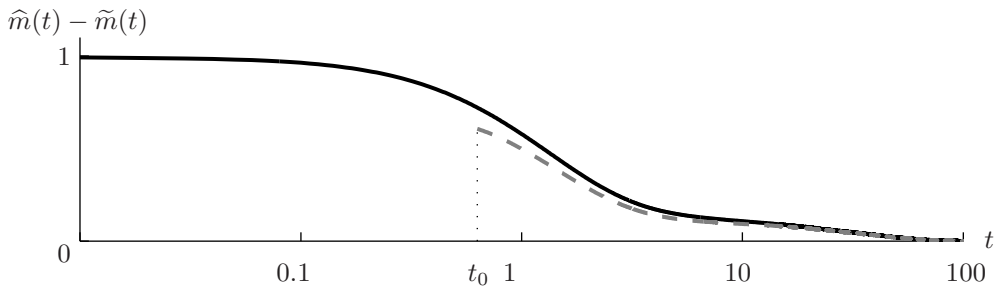
	\caption{Difference between two analytical approximations $\tilde{m}(t)$ and simulation estimate $\hat m(t)$.
The dashed curve uses $D$, while the solid curve does not. The approximations $\tilde{m}(t)$ are easy to evaluate using \eqref{eq:Hunter} together with our proposed approximation \eqref{eq:approx}.
}\label{fig}
\end{figure}

Figure~\ref{fig} plots the difference between the estimated (true value) $\hat m(t)$ and two versions of the approximate $\tilde{m}(\cdot)$.
The solid curve does not use the correction matrix $D$, i.e.\ $\tilde{m}(t)$ is evaluated assuming covariance $C\,t$.
The dashed curve is the improved approximation incorporating an asymptotically exact covariance, $C\,t+D$ (both curves utilize ${\mathbf b}$).
As observed, while both curves converge to zero error as $t \to \infty$, the refinement yields smaller error --- especially in the ``medium'' time horizon after $t_0$.

\section{Proof of Main Result}
\label{sec:proofs}

Without loss of generality we prove Theorem~\ref{thm2} as stated for the case $L=2$. In the ordinary case denote $M(t) \idef \E N(t) $ and the closely related $D_x(t) \idef \E \mathring R_x(t)$ (as well as the $y$-counterpart $D_y(t)$). 

\begin{lemma}
\label{lem:BS4}
For ordinary ${\mathbf R}(\cdot)$:
$
\E \sum_{i < j \leq N(t)-1} X_i\,Y_j = \int_0^t D_{x}(t-s)\,\di D_{y}(s).
$
\end{lemma}
{\bf Proof.} Let $f_n = \di F^{(n)}\,/\,\di M$. Since $M(t) = \sum_{n=0}^\infty F^{(n)}(t)$ for all $t$, $M(t) = 0$ implies all $F^{(n)}(t)=0$, so $F^{(n)} \ll M$ and $f_n$ is well defined. Now $D_{x}(t) =\E \sum_{i=1}^\infty X_i\, \I_{\{S_i\leq t\}} = \sum_{i=1}^\infty \int_0^t  \E\left( X_i\, |\, S_i = s\right)\di F^{(i)}{(s)} = \int_0^t  \sum_{i=1}^\infty \E\left( X_i\, |\, S_i = s\right)f_i(s) \,\di M(s)$, and therefore $\di D_{x}(s)\,/\,\di M(s) = \sum_{i=1}^\infty \E\left(X_i\,|\,S_i=s\right)f_i(s)$. In the following, denote $\tilde S_{j-i} \idef S_j - S_i$. Next,
\begin{align*}
&\E \sum_{i < j \leq N(t)-1} X_i\,Y_j = \sum_{i<j}\E X_i\,Y_j\I_{\{S_j \le t\}} = \sum_{i<j} \E \E\left(X_i\,Y_{j-i}\,\I_{\{\tilde S_{j-i}+S_i \le t\}}\,\Big|\,S_i\right)\\
&=\sum_{i<j} \int_{\omega = 0}^\infty \E\left(X_i\,Y_{j-i}\,\I_{\{\tilde S_{j-i}+\omega \le t\}}\,\Big|\,S_i = \omega\right)\di F^{(i)}(\omega)\\
&=\sum_{i<j}\int_{\omega=0}^t \E\left( X_i\,|\,S_i = \omega\right)\,\E\,\E\left(Y_{j-i}\,\I_{\{\tilde S_{j-i} + \omega \le t\}}\,\Big|\,\tilde S_{j-i}\right)\di F^{(i)}(\omega)\\
&= \sum_{i<j}\int_{\omega=0}^t \E\left( X_i\,|\,S_i = \omega\right)\,\int_{s=\omega}^\infty\E\left(Y_{j-i}\,\I_{\{s \le t\}}\,\Big|\,\tilde S_{j-i}=s-\omega\right)\di F^{(j-i)}(s-\omega)\,\di F^{(i)}(\omega)\\
&= \sum_{i<j}\int_{\omega=0}^t\int_{s=\omega}^t \E\left( X_i\,|\,S_i = \omega\right)\,\E\left(Y_{j-i}\,\Big|\,\tilde S_{j-i}=s-\omega\right)\di F^{(j-i)}(s-\omega)\,\di F^{(i)}(\omega)\\
&= \int_{\omega=0}^t\int_{s=\omega}^t \sum_{i=1}^\infty\sum_{(j-i)=1}^\infty\E\left( X_i\,|\,S_i = \omega\right)\,\E\left(Y_{j-i}\,\Big|\,\tilde S_{j-i}=s-\omega\right)\di F^{(j-i)}(s-\omega)\,\di F^{(i)}(\omega)\\
&= \int_{\omega=0}^t\int_{s=\omega}^t \sum_{i=1}^\infty\E\left( X_i\,|\,S_i = \omega\right)\,\sum_{k=1}^\infty\E\left(Y_{k}\,\Big|\,\tilde S_{k}=s-\omega\right)\,f_{k}(s-\omega)\,\di M(s-\omega)\,f_i(\omega)\,\di M(\omega)\\
&= \int_{\omega=0}^t\int_{s=\omega}^t \sum_{i=1}^\infty\E\left( X_i\,|\,S_i = \omega\right)\,f_i(\omega)\,\sum_{k=1}^\infty\E\left(Y_{k}\,\Big|\,\tilde S_{k}=s-\omega\right)\,f_{k}(s-\omega)\,\di M(s-\omega)\,\di M(\omega)\\
&= \int_{\omega=0}^t D_y(t-\omega)\,\sum_{i=1}^\infty\E\left( X_i\,|\,S_i = \omega\right)\,f_i(\omega)\,\di M(\omega) = \int_{\omega=0}^t D_y(t-\omega)\,\frac{\di\,D_x(\omega)}{\di\,M(\omega)}\,\di M(\omega)\\
&= \int_{0}^t D_y(t-\omega)\,\di D_x(\omega) = \int_{0}^t D_x(t-\omega)\,\di D_y(\omega)\,.
\end{align*}
\begin{raggedright}\halmos\end{raggedright}\\
The next result from \cite{brown1975second} deals with $r_x(t) \idef D_x(t)-a_x\,t-\mathring{b}_x~$:
\begin{lemma}[Restatement of  \cite{brown1975second}, Lemma 3]
\label{lem:BS3}For the ordinary case, if $F(\cdot)$ is spread out and $\mu_3,\,\lambda_1$, and $m_{2,1}$ are finite, then
\[
\int_0^\infty r_x(t) \, \di t = \ell_x\,,
\]
where $\ell_x$ is defined in \eqref{eq:3543}. Moreover, $r_x(\cdot)$ is directly Riemann integrable and $\lim_{t\to\infty} t\,r_x(t) = 0\,$.
\end{lemma}
We can now prove Theorem~\ref{thm2}, which is the key to our approximation \eqref{eq:approx}.

\noindent{\bf Proof of Theorem~\ref{thm2}.}
Since $\E T_1\,|X_1| \leq (\E T_1X_1^2\, \E X_1)^{1/2}$, $\E T_1^2\,|X_1| \leq (\E T_1X_1^2\, \E X_1^3)^{1/2}$, and $\E |X_1Y_1| \leq (\E X_1^2\, \E X_1^2)^{1/2}$, it holds that $m_{1,1}$, $m_{2,1}$, and $p_{0,1,1}$ are finite. Similarly, $n_{1,1}$ and $n_{2,1}$ are finite. It holds,
\begin{align}
\nonumber
\Cov\big(\mathring R_x(t),\,\mathring R_y(t)\big) &=\E \mathring R_{xy}(t) + 2\,\E \sum_{i < j \leq N(t)-1} X_i\,Y_j - \E \mathring R_{x}(t)\,\E \mathring R_{y}(t)\\
\label{eq:1942}
&= \E \mathring R_{xy}(t) +  2\, \int_0^t D_{x}(t-s)\,\di D_{y}(s) - \E \mathring R_{x}(t)\,\E \mathring R_{y}(t)\,,
\end{align}
where the second step follows from Lemma~\ref{lem:BS4}. We have $\E \mathring R_{xy}(t) = a_{xy}\,t + \mathring{b}_{xy} + \lo(1)$ from \eqref{eq:bxyring} and it follows from Theorem~\ref{thm1} that
\begin{align}
\label{eq:Dtilde}
&\E \mathring R_{x}(t)\,\E \mathring R_{y}(t) = a_x\,a_y\,t^2+(a_x\,\mathring b_y+a_y\,\mathring b_x)\,t+\mathring b_x\,\mathring b_y+\lo(1)\,.
\end{align}
Now,
\[
\int_0^t D_{x}(t-s)\,\di D_{y}(s) =\int_0^t r_{x}(t-s)\,\di D_{y}(s)+\int_0^t\big(a_{x}\,(t-s)+\mathring b_{x}\big)\,\di D_{y}(s)\,.
\]
By Lemma~\ref{lem:BS3}, $r_x(\cdot)$ is directly Riemann integrable.  It thus follows from a generalisation of the key renewal theorem to renewal-reward processes (see \cite{brown1972asymptotic}) that $\int_0^t r_{x}(t-s)\,\di D_{y}(s) = a_y\,\ell_x+\lo(1)$.
 Next,
\[
\int_0^t\big(a_x\,(t-s)+\mathring b_x\big)\,\di D_{y}(s) = a_x\,\int_0^t D_{y}(s)\,\di s + \mathring b_x\,D_{y}(t)
=
\mathring b_x\,D_{y}(t)+a_x\,\int_0^t \big(r_y(s) +a_y\,s+\mathring b_y \big)\,\di s\,.\\
\]
Now using Lemma~\ref{lem:BS3} we have
\begin{align}
\label{eq:Dxy2}
\int_0^t D_{x}(t-s)\,\di D_{y}(s) = a_x\,a_y\,t^2\,/2+(a_x\,\mathring b_y+a_y\,\mathring b_x)\,t+a_x\,\ell_y+a_y\,\ell_x+\mathring b_x\,\mathring b_y+\lo(1)\,.
\end{align}
Combining the above into \eqref{eq:1942} yields the result for the ordinary case.\\[1pt]

We now move onto the delayed case.
Since $R_x(t) = \I_{\{T_0\leq t\}}\big(X_0+\mathring{R}_{x}(t-T_0)\big)$ and similarly for $R_y(t)$,
\[
R_x(t)\,R_y(t) = \I_{\{T_0\leq t\}}\big(X_0\,Y_0+X_0\,\mathring{R}_y(t-T_0)+Y_0\,\mathring{R}_x(t-T_0)+\mathring{R}_x(t-T_0)\,\mathring{R}_y(t-T_0)\big).
\]

Now, $\E \I_{\{T_0\leq t\}}X_0Y_0 = \int_0^t \E\left[X_0\,Y_0\,|\,T_0=s\right]\,\di F_0(s) = \E X_0\,Y_0+\lo(1)$.
Next,
\begin{align*}
\E \I_{\{T_0\leq t\}}X_0\,\mathring{R}_y(t-T_0) &= \E\,\I_{\{T_0\leq t\}}\,\E\left[X_0\,|\,T_0\right]\,\big(a_y\,(t-T_0)+\mathring b_y+r_y(t-T_0)\big)\\
&= (a_y\,t+\mathring b_y)\,\E X_0-a_y\,\E T_0\,X_0+\lo(1)\,,
\end{align*}
since $r_y(t)$ converges to 0 as $t\to\infty$ (Theorem~\ref{thm1}) and    both $\sup_t\{|r_y(t)|\}$ and $\E|X_0|$ are finite,  it holds that
$\int_0^t r_y(t-s)\,\E  \left[X_0\,|\,T_0=s\right]\,\di F_0(s) \to 0$ as $t\to\infty$.
Similarly for $\E \I_{\{T_0\leq t\}}Y_0\,\mathring{R}_x(t-T_0)$.

Set $\bar{r}(t) \idef \Cov\big(\mathring{R}_x(t),\,\mathring{R}_y(t)\big)-c_{x,y}\,t-\mathring{d}_{x,y}$. Hence,
\begin{align*}
&\E\,\I_{\{T_0\leq t\}}\mathring{R}_x(t-T_0)\,\mathring{R}_y(t-T_0)\\
&= \E \,\I_{\{T_0\leq t\}}\big(c_{x,y}\,(t-T_0)+d_{x,y}+\bar{r}(t-T_0)\\
&\hspace{55mm}+(a_x\,(t-T_0)+\mathring b_x+r_x(t-T_0))(a_y\,(t-T_0)+\mathring b_y+r_y(t-T_0))\big)\\
&=  \E \,\I_{\{T_0\leq t\}}\big(a_x\,a_y\,t^2+t\,(c_{x,y}-2\,a_x\,a_y\,T_0+\mathring b_y\,a_x+\mathring b_x\,a_y)+\mathring{d}_{x,y}-c_{x,y}T_0+a_x\,a_y\,T_0^2\\
&\hspace{50mm}+\mathring b_x\,\mathring b_y-\mathring b_x\,a_y\,T_0-\mathring b_y\,a_x\,T_0+\bar{r}(t-T_0)+r_x(t-T_0)\,r_y(t-T_0)\\
&\hspace{59mm}+r_x(t-T_0)(a_y\,(t-T_0)+\mathring b_y)+r_y(t-T_0)\,(a_x\,(t-T_0)+\mathring b_x)\big)\,.
\end{align*}
By Theorem~\ref{thm1} and the result proved above for the ordinary case, $r_x(t),\, r_y(t),$ and $\bar{r}(t)$ all converge to~$0$. Moreover, $\sup_t\{|\bar{r}(t)|\},\, \sup_t\{|r_x(t)|\}$, and $\sup_t\{|r_x(t)|\}$ are finite, thus $\int_0^t \bar{r}(t-x)\,\di F_0(x)$, $\int_0^t r_x(t-x)\,r_y(t-x)\,\di F_0(x)$, $\int_0^t r_x(t-x)\,\di F_0(x)$, and $\int_0^t r_y(t-x)\,\di F_0(x)$ all converge to 0 as $t\to\infty$. Further, by Lemma~\ref{lem:BS3}, $t\,r_x(t)$ and  $t\,r_y(t)$ also converge to 0, and it easily follows that $\sup_t\{|t\,r_y(t)|\}$ and $\sup_t\{|t\,r_y(t)|\}$ are finite; thus
\begin{align*}
&\E\,\I_{\{T_0\leq t\}}\big(a_x\,(t-T_0)\,r_y(t-T_0)+a_y\,(t-T_0)\,r_x(t-T_0)\big)\\
 &= \int_0^t \big(a_x\,(t-T_0)\,r_y(t-T_0)+a_y\,(t-T_0)\,r_x(t-T_0)\big)\,\di F_0(x)\,,
\end{align*}
which converges to 0 as $t\to\infty$. Therefore,
\begin{align*}
\E\,\I_{\{T_0\leq t\}}\,\mathring{R}_x(t-T_0)\,\mathring R_y(t-T_0) &= a_x\,a_y\,t^2+(c_{x,y}-2\,a_x\,a_y\,\E T_0+\mathring b_y\,a_x+\mathring b_x\,a_y)\,t+\mathring{d}_{x,y}-c_{x,y}\,\E T_0\\
&\hspace{25mm}+a_x\,a_y\,\E T_0^2+\mathring b_x\,\mathring b_y-\mathring b_x\,a_y\,\E T_0-\mathring b_y\,a_x\,\E T_0+\lo(1)\,.
\end{align*}
Thus,
\begin{align*}
&\E R_x(t)\,R_y(t)\\
&= a_x\,a_yt^2+(c_{x,y}-2\,a_x\,a_y\,\E T_0+\mathring b_y\,a_x+\mathring b_x\,a_y+a_x\,\E Y_0+a_y\,\E X_0)\,t+\mathring{d}_{x,y}-c_{x,y}\,\E T_0+a_x\,a_y\,\E T_0^2\\
&\hspace{5mm}-a_x\,\E T_0\,Y_0-a_y\,\E T_0\,X_0+\E X_0\,Y_0+\mathring b_x\,\mathring b_y-\mathring b_x\,a_y\,\E T_0-\mathring b_y\,a_x\,\E T_0+\mathring b_x\,\E Y_0+\mathring b_y\,\E X_0+\lo(1)\,.
\end{align*}
By Theorem~\ref{thm1},
\begin{align*}
&\E R_x(t)\,\E R_y(t) = (a_x\,t+\mathring b_x+\E X_0-a_x\,\E T_0)\,(a_y\,t+\mathring b_y+\E Y_0-a_y\,\E T_0)+\lo(1)\\
&= a_x\,a_y\,t^2+(a_x\,\mathring b_y+a_x\,\E Y_0-2\,a_x\,a_y\,\E T_0+\mathring b_x\,a_y+a_y\,\E X_0)\,t+\mathring b_x\,\mathring b_y+\mathring b_x\,\E Y_0-\mathring b_x\,a_y\,\E T_0\\
&\hspace{22mm}+\mathring b_y\,\E X_0+\E X_0\,\E Y_0-a_y\,\E X_0\,\E T_0-a_x\,\mathring b_y\,\E T_0-a_x\,\E Y_0\,\E T_0+a_x\,a_y\,(\E T_0)^2+\lo(1)\,.
\end{align*}
Combining the two expressions above yields the result.
\begin{raggedright}\halmos\end{raggedright}

\vspace{-7mm}
\section{Outlook}
\label{sec:conclusion}

The renewal-reward process with multivariate rewards analysed here often plays a role as part of a more complicated stochastic model --- for example in multidimensional risk models, such as in \cite{badila2014two}. Our results may help analysis of such risk models, at least in some asymptotic regime.

We have allowed the distribution of ${\mathbf Z}_n$ to depend on $n$ in a simple way by allowing ${\mathbf Z}_0$ to follow a different distribution to $\{{\mathbf Z}_n\}_{n=1}^\infty$. A possible extension of our work is to allow for more general dependencies of ${\mathbf Z}_n$ on $n$ by partitioning ${\mathbb N}$ into possibly infinite subsets. In \cite{Spataru2010} Sp\u{a}taru gives a CLT for the case of univariate rewards ($L=1$) with unit rewards (a renewal process) in this setting. Extending Sp\u{a}taru's result to the univariate or even the multivariate renewal-reward case remains a challenge.


{\footnotesize \section*{Acknowledgements}
This work was in part carried out as a component of the M.Sc.\ of BP.  YN is supported by Australian Research Council (ARC) grants DP130100156 and DE130100291. }

{\footnotesize
\bibliography{AsymptoticCovarianceRenewalReward}}
\end{document}